\documentstyle[a4,leqno,amsfonts,12pt]{article}
\newcommand{\C}{{\bf C}}

\newcommand{\OC}{\overline{{\bf C}}}

\newcommand{\N}{{\bf N}}
\newcommand{\T}{{\bf T}}

\newcommand{\bP}{{\bf P}}
\newcommand{\D}{{\bf D}}

\newcommand{\de}{\delta}

\newcommand{\mb}{\mbox}

\newcommand{\beq}{\begin{equation}}
\newcommand{\eeq}{\end{equation}}
\newcommand{\oge}{\succeq}
\newcommand{\ole}{\preceq}
\newcommand{\ve}{\varepsilon}

\newcommand{\ov}{\overline}
\newcommand{\al}{\alpha}
\newcommand{\be}{\beta}

\newcommand{\Om}{\Omega}

\newcommand{\om}{\omega}
\newcommand{\z}{\zeta}

\newcommand{\la}{\lambda}
\newcommand{\Ga}{\Gamma}

\newtheorem{th}{Theorem}

\newtheorem{lem}{Lemma}

\newcommand{\ueberschrift}{\bigskip\goodbreak\noindent\bigskip}
\newcounter{theabsatz}
\newcommand{\absatz}[1]{\stepcounter{theabsatz} \ueberschrift
               {\large \bf \arabic{theabsatz}. {#1}} \setcounter{equation}{0}}

\begin{document}

\begin{center}
{\large \bf
Weighted Remez- and Nikolskii-Type 
Inequalities\\ on a Quasismooth Curve
}
\\[2ex] {\bf Vladimir Andrievskii}

\end{center}

\begin{abstract}

We establish sharp $L_p,1\le p<\infty$  weighted Remez- and
Nikolskii-type  inequalities
 for algebraic polynomials considered on a quasismooth (in the sense
 of Lavrentiev) curve in the complex plane.

 \end{abstract}


\footnotetext{
 Received

 Published online
 }

\absatz{Introduction}

From the numerous generalizations of the classical Remez
inequality (see, for example, \cite{rem, borerd,erd,gan}), we mention 
three results which are the starting point of our analysis.

Let $|S|$ be the linear measure (length) of a Borel set $S$ in the
complex plane $\C $. By $\bP_n$ we denote the set of all complex
polynomials of degree at most $n\in\N:=\{ 1,2,\ldots\}$. 
The first
result 
is  due  to Erd\'elyi \cite{erd1}. 
Assume that for $p_n\in\bP_n$ and $\T:=\{ z:|z|=1\}$ we have
\beq \label{1.1} |
\{z\in\T:|p_n(z)|>1\}|
\le s,\quad 0<s\le\frac{\pi}{2} .
 \eeq 
Then,
$|p_n(e^{it})|^2$ is a trigonometric polynomial of degree at
most $n$ and, by the Remez-type inequality on the size of
trigonometric polynomials (cf. \cite[Theorem 2]{erd1} or \cite[p.
230]{borerd}), we obtain 
\beq \label{1.2} 
||p_n||_{C(\T)}\leq e^{2sn},
\quad 0<s\leq\frac{\pi}{2}. 
\eeq
Here $||\cdot||_{C(S)}$ means the uniform norm over $S\subset\C$.

The second result is due to Mastroianni and  Totik \cite{mastot}.
Let $T_n$ be a trigonometric polynomial of degree $n\in \N$, 
$1\le p<\infty$, and $W:[0,2\pi]\to \{x\ge 0\}$ be an $A_\infty$ weight 
function. Then, according to \cite[(5.2) and Theorem 5.2]{mastot}, there are
positive constants $c_1$
and $c_2$
depending only on the $A_\infty$ constant of $W$ and $p$, such that
for a measurable set $E\subset[0,2\pi]$ with $|E|\le s, 0< s\le 1$,
we have
\beq\label{1.3}
\int_{[0,2\pi]}|T_n|^pW\le c_1\exp(c_2sn)
\int_{[0,2\pi]\setminus E}|T_n|^pW.
\eeq
The third result, which is due to Andrievskii and Ruscheweyh \cite{andrus}, extends 
(\ref{1.1})-(\ref{1.2}) to the
case of algebraic polynomials considered on a Jordan curve $\Ga\subset \C$ instead of the
unit circle $\T$.
In the present paper, we always assume that $\Ga$ is  {\it quasismooth} (in the sense of
Lavrentiev), see \cite{pom}, i.e.,
for every $z_1,z_2\in \Ga$,
\beq\label{1.4}
|\Ga(z_1,z_2)|\le \Lambda_\Ga|z_1-z_2|,
\eeq
where $\Ga(z_1,z_2)$ is the shorter arc of $\Ga$ between $z_1$ and $z_2$ (including the  endpoints)
and
$\Lambda_\Ga\ge 1$ is a constant.
 
Let $\Om$ be the unbounded component of $\OC\setminus \Ga$, where $\OC:=\C\cup\{
\infty\}$. Denote by $\Phi$ the conformal mapping of $\Om$ onto
$\D^*:=\{z:|z|>1\}$
with the normalization
$$
\Phi(\infty)=\infty,\quad \Phi'(\infty):=\lim_{z\to\infty}\frac{\Phi(z)}{z}>0.
$$
For $\de>0$ and $A,B\subset\C$, we set
$$
d(A,B)=\mb{dist}(A,B):=\inf_{z\in A,\z\in B}|z-\z|,
$$
$$
\Ga_\de:=\{\z\in\Om:|\Phi(\z)|=1+\de\}.
$$
Let the function $\de(t)=\de(t,\Ga),t>0$ be defined by
the equation $d(\Ga,\Ga_{\de(t)})=t$ and let 
 diam $S$ be the diameter of a set $S\subset\C$.

According to \cite[Theorem 2]{andrus}, if for $p_n\in\bP_n$,
\beq\label{1.5}
|\{z\in\Ga: |p_n(z)|>1\}|\le s<\frac{1}{2}\mb{ diam }\Ga,
\eeq
then
\beq\label{1.6}
||p_n||_{C(\Ga)}\le \exp(c_3\de(s)n)
\eeq
holds with a positive constant $c_3=c_3(\Ga)$.

Our objective is to provide the  weighted $L_p$ analogue of
(\ref{1.5})-(\ref{1.6}) which extends (\ref{1.3}) to the case
of complex polynomials considered on $\Ga$.
Some of our proofs and constructions are  modifications of arguments 
from \cite{mastot, and12, andjat12, els}. For the sake of
completeness, we describe them in detail. 

We denote by $\alpha,
c,\ve,\al_1, c_1,\ve_1,\ldots$ positive constants
(different in different sections)
 that are either  absolute or they depend on
parameters inessential for the argument; otherwise, such  dependence will be
explicitly stated.
For nonnegative
functions $f$ and $g$ we  use the expression $f\ole g$ (order inequality)
if $f\le cg$. The expression $f\asymp g$ means that $f\ole g$ and
$g\ole f$ simultaneously.

\absatz{Main Results}

We say that  a finite Borel measure $\nu$ supported on $\Ga $ is an {\it  $A_\infty$ measure} 
($\nu\in A_\infty(\Ga)$ for short) if
there exists a constant $\la_{\nu}\ge 1$ such that for any
arc $J\subset \Ga$ and a Borel set $S\subset J$ satisfying
$|J|\le 2|S|$ we have
\beq\label{2.1}
\nu(J)\le \la_{\nu}\nu(S),
\eeq
see for instance \cite{coifef, jerken}.
The measure defined by the arclength on $\Ga$ is automatically the
$A_\infty$ measure.
Another interesting example is the equilibrium measure $\mu_\Ga$ on $\Ga$
(see for example \cite{saftot}).
By virtue of \cite{lav}
 $\mu_\Ga\in A_\infty(\Ga)$.
\begin{th}\label{th1}
Let $\nu\in A_\infty(\Ga)$,
$1\le p<\infty$, and let $E\subset \Ga$
be a Borel set.
Then 
for $p_n\in\bP_n,n\in\N$, we have 
\beq\label{2.00}
\int_\Ga |p_n|^p d\nu\le
 c_1\exp(c_2\de(s)n)\int_{\Ga\setminus E}|p_n|^p d\nu
\eeq
provided that $0<|E|\le s < (\em{ diam }\, \Ga)/12$, where the  constants
$c_1$ and $c_2$ depend only on $\Ga, \la_\nu, p$.
\end{th}
Let $\Ga=\T$. Starting with the trigonometric polynomial
$$
T_n(t)=\sum_{k=0}^n(a_k\sin kt+b_k\cos kt),
$$ 
consider the algebraic polynomial
$$
p_{2n}(z):= z^n
\sum_{k=0}^n\left(\frac{a_k}{2i}\left( z^k-z^{-k}\right)+
\frac{b_k}{2}\left( z^k+z^{-k}\right)\right).
$$
Then  (\ref{2.00}) implies (\ref{1.3}) (up to the upper bound on a parameter $s$).

The sharpness of Theorem \ref{th1} is established by our next theorem.
Let $ds=|dz|, z\in\Ga$  be the  arclength measure on
 $\Ga$.
\begin{th} \label{th2}
Let $0<s<\em{ diam }\,\Ga $ and $1\le p<\infty$. Then
there exist an arc $E_s\subset\Ga$ with $|E_s|=s$ as well as constants
$\ve_1=\ve_1(\Ga)$ and $n_0=n_0(s,\Ga,p)\in\N$ such that for any $n>n_0$,
there is a polynomial $p_{n,s}\in\bP_n$ satisfying
\beq\label{2.2n}
\int_\Ga|p_{n,s}|^pds\ge\exp(\ve_1\de(s)n)\int_{\Ga\setminus E_s}|p_{n,s}|^pds.
\eeq
\end{th}
If, in the definition of the $A_\infty$ measure, we assume that $S$ is also an arc,
then $\nu$ is called a {\it doubling measure}.
In \cite[Section 5]{mastot} one can find an example showing that the weighted
Remez-type inequality may not be true in the case of doubling measures.

A straightforward consequence of Theorem \ref{th1} is the following
Nikolskii-type inequality which partially overlaps with  \cite[Corollary 3.10]{var}
where the analogous inequality is proved in another way.
 For more details on the classical Nikolskii
inequality, its generalizations, and further references see, for example \cite{ibrmam, borerd, erd, mam}.
\begin{th}\label{th3}
Let  $\nu\in A_\infty(\Ga)$ satisfy
 $d\nu=wds,w:\Ga\to\{x\ge0\},$
and let $1\le p<q<\infty$.
Then,
for $p_n\in\bP_n,n>n_1$,
\beq\label{2.1s}
\left(\int_\Ga|p_n|^qwds\right)^{1/q}\le c_3 d(\Ga,\Ga_{1/n})^{1/q-1/p}
\left(\int_\Ga|p_n|^p w^{p/q}ds\right)^{1/p}
\eeq
holds with  constants $c_3=c_3(\Ga,p,q,\la_\nu)$ and $n_1=n_1(\Ga)$.
\end{th}
For $\Ga=\T$ (\ref{2.1s}) yields \cite[Theorem 5.5]{mastot}.
The estimate (\ref{2.1s}) is sharp in the following sense.
\begin{th}\label{th4}
For $n\in\N$, there exists a polynomial $p_n^*\in\bP_n$,
such that for $1\le p<q<\infty$,
\beq\label{2.1v}
\left(\int_\Ga|p_n^*|^q ds\right)^{1/q}\ge \ve_2 d(\Ga,\Ga_{1/n})^{1/q-1/p}
\left(\int_\Ga|p_n^*|^p ds\right)^{1/p}
\eeq
holds with $\ve_2=\ve_2(\Ga,p,q).$
\end{th}
Note that $\de(s)$ and $d(\Ga,\Ga_\de)$ can be further estimated.
We mention three well known results. For a more complete theory see,
for example \cite{war, pom1, les, pom}.

The Ahlfors criterion \cite[p. 100]{lehvir} implies that $\Ga$ is quasiconformal.
Therefore, $\Phi$ can be extended to a quasiconformal homeomorphism 
$\Phi:\OC\to\OC$.
Taking into account
Lemma \ref{lem3.1} below and distortion properties of conformal mappings with
quasiconformal extension (cf. \cite[pp. 289, 347]{pom1})
we have 
$$
\delta (s)\ole s^{1/\alpha },\quad 0<s<\mb{diam }\Ga,
$$ 
$$
d(\Ga,\Ga_\de)\oge \de^{\al},\quad 0<\de<1,
$$
with  some $\alpha =\alpha (\Ga)$ such that $1\le\al<2$.

Next, following \cite{pom} we call $\Ga$ {\it Dini-smooth} if
it is smooth and if the angle $\beta(s)$ of the tangent, considered as a
function of the arc length $s$, has the property
   $$
   |\beta(s_2) - \beta(s_1)| \le h(s_2-s_1), \quad 0 < s_2 - s_1 < |\Ga|/2,
   $$
where $h$ is a function satisfying 
   $$
   \int^{|\Ga|/2}_0  \frac{h(x)}{x}\, dx < \infty .
   $$
We call a Jordan arc {\it Dini-smooth} if it is a subarc of some Dini-smooth curve.
According to \cite[Theorem 4]{andrus} if  $\Ga$ is Dini-smooth, then 
$$
\de(s)\asymp s, \quad 0<s<\mb{ diam }\Ga,
$$
$$
d(\Ga,\Ga_\de)\asymp \de ,\quad 0<\de<1.
$$
Moreover, the distortion properties of $\Phi$ in the case of a piecewise
Dini-smooth $\Ga$ (cf. \cite[Chapter 3]{pom} or \cite[pp. 32-36]{andbla}) imply
 that if $\Ga$ consists of a finite number of Dini-smooth arcs 
which meet under the angles $\al_1\pi,\ldots,\al_m\pi$ with respect to $\Om$,
where $0<\al_j<2$, then
$$
\de(s)\asymp s^{1/\al}, \quad 0<s\le\mb{ diam }\Ga,
$$
$$
d(\Ga,\Ga_\de)\asymp \de^{\al},\quad 0<\de<1,
$$
hold with  $\al:=\max(1,\al_1,\ldots,\al_m)$.

\absatz{Auxiliary Constructions and Results}

In this section, we review some of the properties 
of conformal mappings $\Phi$ and $\Psi:=\Phi^{-1}$ whose proofs can be 
found, for example, in \cite[Section 3]{and12}.
We also prove some new facts about these conformal mappings which are used
in the proofs of the main results.

\begin{lem}
\label{lem3.1} 
 Assume  that $z_j\in\overline{\Om},\, t_j:=\Phi(z_j),\,
j=1,2,3$. Then,

(i) the conditions $|z_1-z_2|\ole |z_1-z_3|$ and $|t_1-t_2|\ole
|t_1-t_3|$ are equivalent;

(ii) if $|z_1-z_2|\ole |z_1-z_3|$, then
$$\left|\frac{t_1-t_3}{t_1-t_2}\right|^{1/\al}\ole
\left|\frac{z_1-z_3}{z_1-z_2}\right|\ole
\left|\frac{t_1-t_3}{t_1-t_2}\right|^{\al},
\quad \al=\al(\Ga)\ge1.
$$ 
\end{lem}
Most of the geometrical facts below 
can be obtained by a straightforward application of Lemma \ref{lem3.1} to specifically
chosen triplets of points.

For $\de>0$ and $z\in \Ga$, set
$$
\rho_\de(z):=d(\{z\},\Ga_\de),\quad
\tilde{z}_\de:=\Psi[(1+\de)\Phi(z)].$$
 Then
\beq\label{3.99n} \rho_\de(z)\asymp|z-\tilde{z}_\de|.
\eeq
 Moreover,
for $0<v<u\le 1$ and $z\in \Ga$,
 Lemma \ref{lem3.1} for the triplet
 $z,\tilde{z}_v,\tilde{z}_u$ and (\ref{3.99n}) yield
 \beq\label{3.9n}
\left(\frac{u}{v}\right)^{1/\al}\ole\frac{\rho_u(z)}{\rho_v(z)}\ole
\left(\frac{u}{v}\right)^\al\eeq
which implies
\beq\label{3.2n}
\de(2s)\ole \de(s),\quad 0<s<\mb{ diam }\Ga.
\eeq 
Indeed, the only nontrivial case is where $s$ satisfies $\de(2s)\le1$.
Let $z_{2s}\in\Ga$ be such that
$$
\rho_{\de(2s)}(z_{2s})=d(\Ga,\Ga_{\de(2s)})=2s.
$$
Since
$$
\frac{\rho_{\de(2s)}(z_{2s})}{\rho_{\de(s)}(z_{2s})}\le
\frac{\rho_{\de(2s)}(z_{2s})}{d(L,L_{\de(s)})}
=\frac{2s}{s}=2,
$$
by the left-hand side of (\ref{3.9n}) we obtain  (\ref{3.2n}).

Furthermore, for $0<\de\le1$ and $z,\z\in L$, the following relations hold:

if $|z-\z|\le\rho_\de(z)$, then
 \beq\label{3.8n}
 \rho_\de(\z)\asymp\rho_\de(z);
 \eeq

 if $|z-\z|\ge\rho_\de(z)$, then
 \beq\label{3.1s}
\left(\frac{\rho_\de(z)}{|z-\z|}\right)^{\al_1}\ole
 \frac{\rho_\de(\z)}{|z-\z|}\ole\left(\frac{\rho_\de(z)}{|z-\z|}\right)^{1/\al_1}.
 \eeq
Let $\de_0=\de_0(\Ga)>0$ be fixed such that
 $$
\max_{z\in \Ga}\rho_{\de_0}(z)
 < \frac{|\Ga|}{2}.
 $$
 For $z\in \Ga$ and  $0<\de <\de_0$,
  denote by $z'_\de,z''_\de\in \Ga$ the two points with the properties
 $$
 z\in \Ga(z'_\de,z''_\de),\quad |\Ga (z'_\de,z)|=|\Ga (z,z''_\de)|=\frac{\rho_\de(z)}{2}.
 $$
 If $\de\ge \de_0$, we set $\Ga (z'_\de,z''_\de):=\Ga$. Let
 $$
 l_n(z):=\Ga(z'_{1/n},z''_{1/n}),\quad z\in \Ga,n\in\N.
$$
Hence,  we have
\beq\label{3.1p}
|l_n(z)|\asymp\rho_{1/n}(z).
\eeq
Let $\nu\in A_\infty(\Ga)$.
Consider the function
\beq\label{3.4n}
w_n(z):=\frac{\nu(l_n(z))}{\rho_{1/n}(z)},\quad
z\in \Ga,n\in\N.
\eeq
 Since $\nu$ is also a doubling measure on $\Ga$,
 for any arcs $J_1$ and $J_2$ with
$J_1\subset J_2\subset \Ga$,
\beq\label{3.11n}
\frac{\nu(J_2)}{\nu(J_1)}\le c_1\left(\frac{|J_2|}{|J_1|}\right)^{\al_2},
\quad c_1=c_1(\Ga,\la_\nu), \al_2=\al_2(\Ga,\la_\nu).
\eeq
The proof of (\ref{3.11n}) follows along the same lines as the proof of 
\cite[(4.1)]{andjat12} (cf. \cite[Lemma 2.1]{mastot}).

Next, according to \cite[Lemma 4]{and12}
for $z,\z\in \Ga$ and $n\in\N$,
 \beq\label{3.1k}
\frac{1}{c_2}\left(1+\frac{|\z-z|}{\rho_{1/n}(z)}\right)^{-\al_3}\le
 \frac{w_n(\z)}{w_n(z)}\le c_2\left(1+\frac{|\z-z|}{\rho_{1/n}(z)}\right)^{\al_3},
\eeq
where $c_2=c_2(\Ga,\la_\nu), \al_3=\al_3(\Ga,\la_\nu)$.

We follow a technique of \cite[(3.12)]{and12} and consider for $n,m\in\N$ and $z,\z\in \Ga$
 the polynomial (in $z$)
$$
q_{n,m}(\z,z)=\sum_{j=0}^N a_j(\z)z^j,\quad N=(10n-11)m,
$$
which satisfies the following properties:

 if $|\z-z|\le\rho_{1/n}(z)\asymp \rho_{1/n}(\z)$, then
 \beq\label{3.13n}
 \frac{1}{c_3}\le|q_{n,m}(\z,z)|\le c_3,\quad c_3=c_3(\Ga,m);
 \eeq

 if $|\z-z|>\rho_{1/n}(z)$, then by virtue of (\ref{3.1s}),
 \beq\label{3.14n}
 |q_{n,m}(\z,z)|
 \le c_4
\left(
 \frac{\rho_{1/n}(\z)}{|\z-z|}\right)^{ m}
\le c_5
 \left(
 \frac{\rho_{1/n}(z)}{|\z-z|}\right)^{m/\al_1},
\eeq
where
 $c_j=c_j(\Ga,m), j=4,5$.

 Let  for  $z\in\Ga $ and $n,m\in\N$,
\begin{eqnarray}
I_{n,m}(z)&:=&\int_\Ga|q_{n,m}(\z,z)|
\frac{w_n(\z)}{w_n(z)}
\frac{|d\z|}{\rho_{1/n}(\z)}\nonumber\\
&=&\label{3.15n}
\frac{1}{w_n(z)}\int_\Ga |q_{n,m}(\cdot,z)|\frac{w_n}{\rho_{1/n}} ds.
\end{eqnarray}
We use the following notation: for $z\in\C$ and $\de>0$,
$$
D(z,\de):=\{\z:|\z-z|<\de\},\quad 
D^*(z,\de):=\C\setminus\ov{D(z,\de)}.
$$
\begin{lem}\label{lem3.2}
There exist sufficiently large $m=m(\Ga,\la_\nu)\in\N$ and $c_6=c_6(\Ga,\la_\nu)$ such that
\beq\label{3.16n}
\frac{1}{c_6}\le I_{n,m}(z)
\le c_6,\quad z\in \Ga.
\eeq
\end{lem}
{\bf Proof}. 
According to the inequalities (\ref{1.4}),
 (\ref{3.8n}), (\ref{3.1k}), and (\ref{3.13n})  we obtain
$$
I_{n,m}(z)
\ge
\int_{\Ga\cap D(z,\rho_{1/n}(z))}|q_{n,m}(\z,z)|
\frac{w_n(\z)}{w_n(z)}
\frac{|d\z|}{\rho_{1/n}(\z)}
\oge 1,
$$
which yields the left-hand side of (\ref{3.16n}).

Next, by (\ref{1.4}), (\ref{3.8n}), (\ref{3.1s}), (\ref{3.1k})-(\ref{3.14n}),
and \cite[(3.20)]{and12}
we have
\begin{eqnarray*}
I_{n,m}(z)
&\ole&
\int_{\Ga\cap \ov{D(z,\rho_{1/n}(z)})}|q_{n,m}(\z,z)|
\frac{w_n(\z)}{w_n(z)}
\frac{|d\z|}{\rho_{1/n}(z)}\\
&&+
\int_{\Ga\cap D^*(z,\rho_{1/n}(z))}
|q_{n,m}(\z,z)|
\frac{w_n(\z)}{w_n(z)}\frac{|\z-z|}{\rho_{1/n}(\z)}
\frac{|d\z|}{|\z-z|}
\\
&\ole& 1+\int_{\Ga\cap D^*(z,\rho_{1/n}(z))}
 \left(\frac{|\z-z|}{\rho_{1/n}(z)}\right)^{\al_3-m/\al_1+\al_1}\frac{|d\z|}
{|\z-z|}\ole 1
 \end{eqnarray*}
if $m$ is any (fixed) number with $\al_3-m/\al_1+\al_1<0$.

Hence, the right-hand side of (\ref{3.16n}) is also proved.

\hfill$\Box$

\begin{lem}\label{lem3.1n}
For $r\ge1$,
\beq\label{3.1v}
\frac{\rho_{1/n}(\z)}{c_7}\le\int_\Ga|q_{n,2}(\z,z)|^r|dz|
\le c_7\rho_{1/n}(\z),\quad \z\in\Ga,
\eeq
where $c_7=c_7(\Ga,r)$.
\end{lem}
{\bf Proof}.
The left-hand side inequality follows from (\ref{1.4}) and (\ref{3.13n}):
$$
\int_\Ga|q_{n,2}(\z,z)|^r|dz|
\ge
\int_{\Ga\cap D(\z,\rho_{1/n}(\z))}|q_{n,2}(\z,z)|^r|dz|
\oge\rho_{1/n}(\z).
$$
Furthermore, according to (\ref{1.4}), (\ref{3.13n}), (\ref{3.14n}), and \cite[(3.20)]{and12}
we have
\begin{eqnarray*}
&&\int_\Ga|q_{n,2}(\z,z)|^r|dz|\\
&\ole&
\int_{\Ga\cap \ov{D(\z,\rho_{1/n}(\z)})}|dz|
+
\int_{\Ga\cap D^*(\z,\rho_{1/n}(\z))}
 \left(\frac{\rho_{1/n}(\z)}{|\z-z|}\right)^{2r}|d\z|\\
&\ole& \rho_{1/n}(\z),
 \end{eqnarray*}
which proves the right-hand side of (\ref{3.1v}).

\hfill$\Box$

\absatz{Proofs of Theorems}

We start with some preliminaries.
Let 
$$
q_r(z):=c\prod_{j=1}^m|z-z_j|^{\be_j},\quad z\in\C,
$$
where $z_j\in\C,c>0,\be_j>0$ be
a {\it generalized polynomial} 
of degree $r:=\be_1+\ldots+\be_m$
and let
$$
E(q_r):=\{z\in \Ga:q_r(z)>1\}.
$$
By \cite[Theorem 2]{andrus}, the condition
\beq\label{4.0t}
|E(q_r)|\le s<\frac{1}{2}\mb{ diam }\Ga
\eeq
yields
\beq\label{4.1}
||q_r||_{C(\Ga)}\le\exp(c_1\de(s)r),\quad c_1=c_1(\Ga).
\eeq
Consider the set
$$
F_s=F_{s}(q_r):=\{z\in \Ga:q_r(z)>\exp(-c_1\de(s)r)||q_r||_{C(\Ga)}\}
$$
and the generalized polynomial
$$
f_{r,s}(z):=\frac{q_r(z)\exp(c_1\de(s)r)}{||q_r||_{C(\Ga)}}
$$
so that $E(f_{r,s})=F_{s}$.

We have
\beq\label{4.2}
|F_s|\ge s,\quad 0<s< \frac{1}{2}\mb{ diam }\Ga.
\eeq
Indeed, the case $|F_s|\ge($diam $\Ga)/2$ is trivial.
If
$|F_{s}|<($ diam $\Ga)/2$, then by (\ref{4.0t})-(\ref{4.1}), applied to $f_{r,s}$
and $|F_s|$ instead of $q_r$ and $s$,
we obtain
$$
\exp(c_1\de(s)r)=||f_{r,s}||_{C(\Ga)}\le\exp
\left(c_1\de(|F_s|)r\right),
$$
that is,
$
\de(|F_s|)\ge\de(s)
$
which  implies (\ref{4.2}).

Let, as before, $1\le p<\infty$. We claim that if a Borel set $A\subset\Ga$
satisfies
\beq\label{4.1w}
|A|\ge|\Ga|-s,\quad 0<s<\frac{1}{4}\mb{ diam }\Ga,
\eeq
then
\beq\label{4.3}
\int_{\Ga}(q_r)^p ds\le \left( 1+\exp(c_2\de(s)pr)\right)\int_{A}(q_r)^p ds,\quad
c_2=c_2(\Ga).
\eeq
Indeed, by virtue of (\ref{4.2}) for $0<s< ($diam $\Ga)/4$, we have
$|F_{2s}|\ge 2s$ which yields $ |A\cap F_{2s}|\ge s.$
Therefore, according to (\ref{3.2n}),
\begin{eqnarray*}
\int_{\Ga\setminus A}(q_r)^p ds&\le&
s||q_r||_{C(\Ga)}^p\le \int_{A\cap F_{2s}}||q_r||^p_{C(\Ga)}ds
\\
&\le&\exp\left(c_1\de\left(2s\right)pr\right)\int_{A\cap F_{2s}}(q_r)^p ds\\
&\le&
\exp\left(c_2\de\left(s\right)pr\right)\int_{A}(q_r)^p ds,
\end{eqnarray*}
which proves (\ref{4.3}).

Let $w_n,n\in\N$ be defined by (\ref{3.4n}).
\begin{lem}\label{lem4.1}
For a Borel set $A\subset \Ga$ satisfying (\ref{4.1w}),
$ 1\le p<\infty$, and
$p_n\in\bP_n,n\in\N$,
\beq\label{4.4}
\int_\Ga|p_n|^pw_n ds\le c_3\exp(c_4\de(s)n)
\int_A|p_n|^pw_n ds,
\eeq
where $c_j=c_j(\Ga,p,\la_\nu),j=3,4.$
\end{lem}
{\bf Proof}.
Let $q_{n,m}$ be the polynomial defined in Section 3.
 By  (\ref{4.3}) applied to the generalized polynomial
 $q_r:=|p_n| |q_{n,m}(\z,\cdot)|^{1/p}$, where $\z\in \Ga, m=m(\Ga)$ is from Lemma
 \ref{lem3.2} and $r\asymp n$, we have
 $$
 \int_\Ga|p_n|^p|q_{n,m}(\z,\cdot)|ds\le
 (1+\exp (c_4\de(s)n)) 
 \int_A|p_n|^p|q_{n,m}(\z,\cdot)| ds.
$$
Multiplying the both sides of this inequality by 
$w_n(\z)/\rho_{1/n}(\z)$,
integrating by $\z$ over $\Ga$, and applying the Fubini theorem we obtain
 \begin{eqnarray*}
 &&\int_\Ga|p_n|^pw_n I_{n,m} ds\\
&=&\int_\Ga\int_\Ga |p_n(z)|^p|q_{n,m}(\z,z)|\frac{w_n(\z)}{\rho_{1/n}(\z)}|d\z||dz|\\
&\le& (1+\exp(c_4\de(s)n)) \int_A\int_\Ga
|p_n(z)|^p|q_{n,m}(\z,z)|\frac{w_n(\z)}{\rho_{1/n}(\z)}|d\z||dz|\\
 &=&
(1+\exp (c_4\de(s)n)) \int_A|p_n|^pw_nI_{n,m}ds,
 \end{eqnarray*}
 where $I_{n,m}$ is defined by (\ref{3.15n}),
which, together with (\ref{3.16n}), yields (\ref{4.4}).

\hfill$\Box$

{\bf Proof of Theorem \ref{th1}}. 
The construction below is partially adapted from the proof of \cite[Theorem 3.1]{mastot}
and the proof of \cite[Theorem 4]{andjat12}.
  Let $m=m(n,s,\Ga)\in\N$ be a sufficiently large number  to be chosen later and
  let 
 $$
 \theta_k:=\frac{2\pi k}{N},\quad \xi_k:=\Psi(e^{i\theta_k}),
\quad k=0,\ldots,N:=nm,
 $$
$$
J_k':=\{e^{i\theta}:\theta_{k-1}\le\theta<\theta_k\},
\quad J_k:=\Psi(J_k'),\quad k=1,\ldots,N.
$$
  By virtue of Lemma \ref{lem3.1},  (\ref{1.4}), and (\ref{3.99n})
  for  $ k=1,\ldots,N,$ we have
  \beq\label{4.5}
  |J_k|\asymp|\xi_k-\xi_{k-1}|\asymp|\xi_k-(\tilde{\xi_k})_{1/N}|
  \asymp\rho_{1/N}(\xi_k).
  \eeq
Let $K:=\{ k:|E\cap J_k|<|J_k|/2\}$. Then
$$
\sum_{k\not\in K}|J_k|\le 2\sum_{k\not\in K}|E\cap J_k|\le 2|E|<2s,
$$
which for 
$$
A:=\Ga\setminus E\quad \mb{and}\quad A^*:=\bigcup_{k\in K}(A\cap J_k)\subset A
$$
implies
\begin{eqnarray*}
|\Ga|&=& \bigcup_{k\not\in K}|A\cap J_k|+|A^*|+|E|\\
&\le& \bigcup_{k\not\in K}|J_k|+|A^*|+|E|<|A^*|+3s,
\end{eqnarray*}
that is,
\beq\label{4.6}
|A^*|>|\Ga|-3s.
\eeq
Let 
$$
B_k:=\sup_{\xi\in J_k}w_n(\xi),
$$
and let $\eta_k,v_k\in \ov{J_k}$ be such that
$$
|p_n(\eta_k)|=\min_{\xi\in \ov{J_k}}|p_n(\xi)|,\quad
|p_n(v_k)|=\max_{\xi\in \ov{J_k}}|p_n(\xi)|=||p_n||_{C(J_k)}.
$$
Consider
$$
R=R(p_n, p, m,n):=\sum_{k\in K}|p_n(v_k)|^pB_k|J_k\cap A|
$$
and
$$
V=V(p_n, p, m,n):=R-\sum_{k\in K}|p_n(\eta_k)|^pB_k|J_k\cap A|
$$
which satisfy
\begin{eqnarray*}
V&=& \sum_{k\in K}(|p_n(v_k)|^p-|p_n(\eta_k)|^p) B_k|J_k\cap A|\\
&\le&
p\sum_{k\in K} |p_n(v_k)-p_n(\eta_k)| |p_n(v_k)|^{p-1}B_k|J_k\cap A|.
\end{eqnarray*}
If $p>1$ and $q>1$ satisfy $1/p+1/q=1$, H\"{o}lder's inequality implies
 \begin{eqnarray*}
 V&\ole&\left(\sum_{k\in K} |p_n(v_k)-p_n(\eta_k)|^p B_k|J_k\cap A|\right)^{1/p}
 R^{1/q}\\
 &\le&
 \left(\sum_{k\in K}\left(\int_{J_k}|p_n'|ds\right)^p
 B_k|J_k\cap A|\right)^{1/p}R^{1/q}.
 \end{eqnarray*}
 If $p=1$, setting  $R^{1/q}:=1$, we have  the same estimate for $V$.

Note that by (\ref{3.1k}) and (\ref{4.5}) $B_k\asymp A_k:=\inf_{\xi\in J_k} w_n(\xi)$.

 Since  H\"{o}lder's inequality also yields
 $$
 \left(\int_{J_k}|p_n'|ds\right)^p\le|J_k|^{p-1}\int_{J_k}|p_n'|^p ds,
 $$
 by \cite[Lemma 1]{and12}, Lemma \ref{lem3.1},  (\ref{3.9n})-(\ref{3.8n}), 
(\ref{3.1k}), and (\ref{4.4})-(\ref{4.6})
 for the nonzero polynomial $p_n$  we further have
 \begin{eqnarray*}
 VR^{-1/q}&\ole& \left(\sum_{k=1}^N\left(\frac{\rho_{1/N}(\xi_k)}{\rho_{1/n}(\xi_k)}\right)^p
 \int_{J_k}(|p_n'|
 \rho_{1/n})^{p}w_n ds\right)^{1/p}\\
 &\ole&
 m^{-\ve}\left(
 \int_{\Ga}(|p_n'|
 \rho_{1/n})^{p}w_n ds\right)^{1/p}\\
 &\ole&
m^{-\ve }\left(
 \int_{\Ga}|p_n|^p
 w_n ds\right)^{1/p}\\
&\ole& m^{-\ve}\exp(c_4\de(3s)n)
\left(\int_{A^*}|p_n|^pw_n ds\right)^{1/p}\\
&\le&
m^{-\ve}\exp(c_5\de(s)n)R^{1/p},
 \end{eqnarray*}
i.e.,
$$
V\le c_6 m^{-\ve}\exp(c_5\de(s)n)R.
$$
Taking $m$ to be  the integral part of
$$
1+\left(2c_6 \exp(c_5\de(s)n)\right)^{1/\ve}
$$
we have $ V\le R/2$ and
$
m\asymp \exp(c_7\de(s)n).
$
Therefore, 
$$
R\le
2\sum_{k\in K}|p_n(\eta_k)|^p B_k|J_k \cap A|
\asymp \sum_{k\in K}|p_n(\eta_k)|^p A_k|J_k \cap A|.
$$
Since $\nu\in A_\infty(\Ga)$ and
$$
|J_k\cap A|\ge \frac{|J_k|}{2},\quad k\in K,
$$
according to (\ref{2.1}),  (\ref{3.1p}), (\ref{3.11n}), and (\ref{4.5})
for $\xi\in J_k$ we have
\begin{eqnarray*}
w_n(\xi)&=&\frac{\nu(l_n(\xi))}{\rho_{1/n}(\xi)}\\
&\ole&
\frac{1}{|J_k|}\frac{\nu(l_n(\xi))}{\nu(l_N(\xi))}
\frac{\nu(l_N(\xi))}{\nu(J_k)}
\frac{\nu(J_k)}{\nu(J_k\cap A)}\nu(J_k\cap A)\\
&\ole& m^{\al_1}\frac{\nu(J_k\cap A)}{|J_k|}.
\end{eqnarray*}
Therefore,
$$
R\ole m^{\al_1}\sum_{k\in K}|p_n(\eta_k)|^p\nu (J_k\cap A)
\le m^{\al_1}\int _{A^*} |p_n|^p d\nu.
$$
Moreover, by \cite[Lemma 2]{and12}, (\ref{3.2n}), (\ref{4.6}), and Lemma \ref{lem4.1},
\begin{eqnarray*}
\int_{\Ga}|p_n|^pd\nu&\ole&\int_\Ga|p_n|^pw_nds\ole \exp(c_4\de(3s)n)\int_{A^*}
|p_n|^pw_nds\\
&\ole&\exp(c_5\de(s)n)R
\ole \exp( (c_5+c_7\al_1)\de(s)n)\int_{A^*}|p_n|^pd\nu\\
&\le& \exp(c_{8}\de(s)n)\int_A|p_n|^pd\nu,\quad c_{8}=c_5+c_7\al_1,
\end{eqnarray*}
which is the desired conclusion.

\hfill$\Box$

{\bf Proof of Theorem \ref{th2}.}
Let $z_s\in\Ga$ and $\z_s\in \Ga_{\de(s)}$
satisfy $|z_s-\z_s|=s=\rho_{\de(s)}(z_s).$
Define points $z_s^*,z_s^{**}\in \Ga$ such that $z_s\in
\Ga(z_s^*,z_s^{**})=:E_s$ and
$$
|\Ga(z_s^*,z_s)|=|\Ga(z_s,z_s^{**})|=\frac{s}{2}\, ,
$$
i.e., $|E_s|=s$.

Lemma \ref{lem3.1} and (\ref{1.4}) yield
\beq\label{4.1p}
|\Phi(z_s^*)-\Phi(z_s)|\asymp |\Phi(z_s)|-\Phi(z_s^{**})|\asymp
|\Phi(z_s) -\Phi(\z_s)|\ge \de(s).
\eeq
Let $A_s:=\ov{\Ga\setminus E_s}$ and let $\Phi_s$ be the conformal mapping of
$\Om_s:=\OC\setminus A_s$ onto $\D^*$ normalized by
$$
\Phi_s(\infty)=\infty,\quad \Phi_s'(\infty)>0.
$$
According to \cite[Lemma 5]{andrus}, (\ref{1.4}),  and (\ref{4.1p}) we obtain
$$
\log|\Phi_s(z_s)|\oge |\Phi_s(z_s^{**})-\Phi_s(z_s^*)|\oge \de(s),
$$
that is,
\beq\label{4.11p}
|\Phi_s(z_s)|\ge \exp(\ve_1\de(s))\ge 1+\ve_1\de(s),\quad
\ve_1=\ve_1(\Ga).
\eeq
Let $p_{n,s}\in\bP_n$ be the $n$-th Faber polynomial for $\Om_s$
(see \cite[Chapter II, \S 1]{smileb} or \cite[Chapter II]{sue}).
From a result by Pommerenke \cite[p. 85, Theorem 3.11]{pom1} (see
also \cite[Chapter IX, \S 3]{sue}) it follows that
\beq\label{4.2p}
||p_{n,s}||_{C(A_s)}\le 2\sqrt{n(\log n+2)}.
\eeq
Moreover, according to \cite[Chapter II, \S 1]{smileb} for $\xi\in\Om_s$,
\beq\label{4.22p}
p_{n,s}(\xi)=\Phi_s(\xi)^n+\om_{n,s}(\xi),
\eeq
where
\beq\label{4.3p}
|\om_{n,s}(\xi)|\le\left(n\log\frac{|\Phi_s(\xi)|^2}{|\Phi_s(\xi)|^2-1}
\right)^{1/2}.
\eeq
Next, by (\ref{1.4}) for $d_s:=$dist$(z_s,A_s)$ we have
$\ve_2 s\le d_s\le s/2$, where $\ve_2=\ve_2(\Ga)$.

According to \cite[p. 23, Lemma 2.3]{andbla} for
$\xi\in W_s:=\Ga\cap D(z_s,d_s/32)$ we obtain
$$
|\Phi_s(\xi)-\Phi_s(z_s)|\le\frac{1}{2}(|\Phi_s(z_s)|-1),
$$
and by (\ref{4.11p})
\beq\label{4.4p}
|\Phi_s(\xi)|\ge 1+\frac{\ve_1}{2}\de(s),
\eeq
which, together with (\ref{4.3p}), implies
\begin{eqnarray}
||\om_{n,s}||_{C(W_s)}&\le&\left( n\log\left(1+\frac{1}{(1+\frac{\ve_1}{2}\de(s))^2-1}
\right)\right)^{1/2}
\nonumber\\
\label{4.5p}
&\le&\left( n\log\left(1+\frac{1}{\ve_1 \de(s)}
\right)\right)^{1/2}.
\end{eqnarray}
Furthermore, (\ref{4.22p}), (\ref{4.4p}), and   (\ref{4.5p}) yield
$$
||p_{n,s}||_{C(W_s)}\ge\left(1+\frac{\ve_1}{2}\de(s)\right)^n
-\left( n\log\left(1+\frac{1}{\ve_1 \de(s)}
\right)\right)^{1/2}.
$$
Let $n_2=n_2(\Ga,s)\in\N$ and $\ve_3=\ve_3(\Ga)$ be such that 
$$
||p_{n,s}||_{C(W_s)}\ge\frac{1}{2}\left(1+\frac{\ve_1}{2}\de(s)\right)^n,\quad
n>n_2,
$$
$$
1+\frac{\ve_1}{2}\de(s)\ge\exp(2\ve_3\de(s)),
$$
that is,
$$
||p_{n,s}||_{C(W_s)}\ge\frac{1}{2}\exp(2\ve_3\de(s)n),\quad
n>n_2.
$$
Summarizing, by virtue of (\ref{4.2p}), we have
\begin{eqnarray*}
&&\frac{\exp(\ve_3\de(s)n)\int_{A_s}|p_{n,s}|^pds}{\int_{W_s}|p_{n,s}|^pds}\\
&\le&
\frac{\exp(\ve_3\de(s)n)|\Ga|2^p(n\log(n+2))^{p/2}}{2^{-p}
\exp(2\ve_3\de(s)np)\ve_216^{-1}s}\\
&=& 4^{p+2}|\Ga| (n\log(n+2))^{p/2} s^{-1}\ve_2^{-1}\exp(-\ve_3\de(s)n)\to 0\quad\mb{as }n\to\infty.
\end{eqnarray*}
Let $n_0=n_0(s,\Ga,p)>n_2$ be such that for $n>n_0$ the right-hand side of the
last inequality is at most $1$. Then, the left-hand side is also $\le1$
from  which
(\ref{2.2n}) follows.

\hfill$\Box$

{\bf Proof of Theorem \ref{th3}}.
Modifying the reasoning from the proof of \cite[Theorem 5.5]{mastot}, we
let $d_n:=d(\Ga,\Ga_{1/n})$ and
$$
E_n=E_{n,q}:=\left\{z\in\Ga:|p_n(z)|^qw(z)\ge d_n^{-1}\int_\Ga|p_n|^qwds\right\}.
$$
Since
$$
\int_\Ga|p_n|^q wds\ge |E_n|d_n^{-1}\int_\Ga|p_n|^q wds
,
$$
we have $|E_n|\le d_n$. 

According to (\ref{1.4}) and Lemma \ref{lem3.1}, there exists $n_1=n_1(\Ga)\in \N$
such that for $n>n_1$ we have $d_n<($diam $\Ga)/12$.

Since $\de(d_n)=1/n$,  by Theorem \ref{th1} for $n>n_1$ we obtain
\begin{eqnarray*}
\int_\Ga|p_n|^q wds&\ole& \int_{\Ga\setminus E_n}|p_n|^qwds
=\int_{\Ga\setminus E_n}(|p_n|^p w^{p/q})(
|p_n|^qw)^{(q-p)/q}ds\\
&\le& 
\left(d_n^{-1}\int_\Ga|p_n|^qwds\right)^{(q-p)/q}
\int_\Ga|p_n|^pw^{p/q}ds,
\end{eqnarray*}
that is,
$$
\left(\int_\Ga|p_n|^qwds\right)^{p/q}\ole
d_n^{(p-q)/q}\int_\Ga|p_n|^pw^{p/q}ds,
$$
which establishes (\ref{2.1s}).

\hfill$\Box$

{\bf Proof of Theorem \ref{th4}}.
There is no loss of generality in assuming that $n>100$
(for $n\le 100$ take $p_n^*\equiv 1$).
Let $z_{1/n}\in\Ga$ satisfy
$$
\rho_{1/n}(z_{1/n})=\min_{z\in\Ga}\rho_{1/n}(z)=
d(\Ga,\Ga_{1/n}).
$$
Consider polynomial 
$$
p_n^*(z):=q_{k,2}(z_{1/n},z),
$$ where $k$ is the integral part of $n/20$ and $q_{k,2}$
is introduced in Section 3.
By (\ref{3.9n}) and Lemma \ref{lem3.1n} for any fixed $r\ge 1$,
$$
\int_\Ga|p_n^*|^rds\asymp\rho_{1/k}(z_{1/n})\asymp
d(\Ga,\Ga_{1/n}),
$$
which implies (\ref{2.1v}).

\hfill$\Box$

{\bf Acknowledgements}
Part of this work was done during the Fall of 2016 semester, while the author visited
the Katholische Universit\"at Eichst\"att-Ingolstadt and
the Julius Maximilian University of W\"urzburg.
 The author is  also grateful to   M.
 Nesterenko
 for his helpful comments.

V. V. Andrievskii

 Department of Mathematical Sciences

 Kent State University

 Kent, OH 44242

 USA

e-mail: andriyev@math.kent.edu

\end{document}